\documentclass[12pt]{amsart}

\def\H{\mathbb{H}}
\def\C{\mathbb{C}}
\def\Z{\mathbb{Z}}
\def\N{\mathbb{N}}

\def\R{\mathbb{R}}

 \newtheorem{thm}{Theorem}[section]
 \newtheorem{cor}[thm]{Corollary}
 \newtheorem{lem}[thm]{Lemma}

\newcommand{\be}{\begin{equation}}
\newcommand{\ee}{\end{equation}}
\newcommand{\bea}{\begin{eqnarray}}

\newcommand{\eea}{\end{eqnarray}}
\newcommand{\Bea}{\begin{eqnarray*}}
\newcommand{\Eea}{\end{eqnarray*}}

\newcounter{cnt1}
\newcounter{cnt2}
\newcounter{cnt3}
\newcommand{\blr}{\begin{list}{$($\roman{cnt1}$)$}
 {\usecounter{cnt1} \setlength{\topsep}{0pt}
 \setlength{\itemsep}{0pt}}}
\newcommand{\bla}{\begin{list}{$($\alph{cnt2}$)$}
 {\usecounter{cnt2} \setlength{\topsep}{0pt}
 \setlength{\itemsep}{0pt}}}
\newcommand{\bln}{\begin{list}{$($\arabic{cnt3}$)$}
 {\usecounter{cnt3} \setlength{\topsep}{0pt}
 \setlength{\itemsep}{0pt}}}
\newcommand{\el}{\end{list}}

\date{}
\begin{document}

\title[ Gutzmer's formula ]
{ An analogue of Gutzmer's formula for  \\
\vskip .5em  Hermite expansions\\
\vskip 1.5em {\tt by} }

\author[Thangavelu]{S.\ Thangavelu}
\address{Department of Mathematics\\
Indian Institute of Science\\
Bangalore 560 012, India. {\it e-mail~:} {\tt
veluma@math.iisc.ernet.in}}

\keywords{ Hermite and Laguerre functions, Heisenberg group, symplectic, orthogonal and unitary matrices, entire functions}
\subjclass{Primary: 42A38 ; Secondary: 42B99, 43A90}

\maketitle

\begin{abstract}
We prove an analogue of Gutzmer's formula for Hermite expansions. As a consequence we obtain a new proof of a characterisation of the image of $ L^2(\R^n)$ under the Hermite semigroup. We also obtain some new orthogonality relations for complexified Hermite functions.
\end{abstract}

\section{\bf Introduction }
\setcounter{equation}{0}

By Gutzmer's formula we mean any analogue of the  formula
$$ (2\pi)^{-1}\int_0^{2\pi} |f(x+iy)|^2 dx = 
\sum_{k=-\infty}^\infty |\hat{f}(k)|^2 
e^{-2ky} $$
valid for any $ 2\pi $ periodic holomorphic function $ f $ in a strip in the
complex plane. Here $ \hat{f}(k) $ stands for the Fourier coefficients of the restriction of $ f $ to the real line. An analogue of such a formula was established by Lassalle [9] for holomorphic functions on the complexification of compact symmetric spaces. A similar formula for holomorphic functions on the complex
crowns associated to noncompact Riemannian symmetric spaces was discovered by
Faraut [3]. As can be seen from Faraut [4] and Kr\"otz-Olafsson-Stanton [8] 
such 
formulas are useful in the study of Segal-Bargmann or heat kernel transforms.
\newpage
Recently in [15] we have proved an analogue of Gutzmer's formula on the
Heisenberg groups and used them to study heat kernel transforms and 
Paley-Wiener theorems.

In this paper we prove an analogue of Gutzmer's formula for Hermite expansions.
Let $ H $ be the Hermite operator on $ \R^n $ having the spectral decomposition
$ H = \sum_{k=0}^\infty (2k+n) P_k.$ Let $ \H^n = \R^n \times \R^n \times \R $
be the Heisenberg group whose complexification is $ \C^n \times \C^n \times \C.$ Let $ \pi(x,u) $ be the projective representation of $ \R^n \times \R^n $
related to the Schr\"odinger representation of $ \H^n $ and denote by 
$ \pi(x+iy,u+iv) $ its extension to $ \C^n \times \C^n.$ Let $ K = 
Sp(n,\R)\cap O(2n,\R) $ which acts on $ \C^n \times \C^n.$ Denote by
 $ \varphi_k(z,w) $ the Laguerre functions
of type $ (n-1) $ extended to $ \C^n \times \C^n.$ Our main result is the following.

\begin{thm} Let $ F $ be an entire function on $ \C^n .$ Denote by $ f $ its restriction to $ \R^n.$ Then for any  $ z = x+iy, w= u+iv \in \C^n $  we have
$$ \int_{\R^n} \int_K |\pi(\sigma.(z,w))F(\xi)|^2 d\sigma d\xi $$
$$ =  e^{(u\cdot y - v \cdot x)} \sum_{k=0}^\infty \frac{k!(n-1)!}{(k+n-1)!}
\varphi_k(2iy,2iv) \|P_kf\|_2^2 .$$
\end{thm}

As an immediate corollary we obtain the following characterisation of the image
of $ L^2(\R^n) $ under the Hermite semigroup $ e^{-tH}, t > 0.$ Let 
$$ U_t(x,y) = 2^n(\sinh(4t))^{-\frac{n}{2}}e^{\tanh(2t)|x|^2 -\coth(2t)|y|^2}.
$$

\begin{cor} An entire function $ F $ on $ \C^n $ belongs to the image of $ L^2(\R^n) $ under $ e^{-tH} $ if and only if
$$ \int_{\R^n}\int_{\R^n} |F(x+iy)|^2 U_t(x,y) dx dy < \infty.$$
\end{cor}

This characterisation is not new and there are several proofs available in the
literature, see Byun [1], Karp [6] and Thangavelu [14]. In Section 4 we derive 
some more consequences of the Gutzmer's fomula. 

We conclude the introduction with some remarks about the methods used in proving Gutzmer formulas. As in the case of Fourier series, Lassalle [9] used Plancherel theorem for the Laurent expansions of  holomorphic functions on the 
complexifications of compact symmetric spaces $ X = K/M.$  The matrix coefficients associated to class one represenations in the unitary dual of a compact Lie group $ K $ holomorphically extend to its
complexification $ K_\C.$ Thus any function $ f $ whose 'Fourier coefficients' have exponential decay can be extended to the 
complexification $ X_\C = K_\C/M_\C.$
Then by appealing to Plancherel theorem and using orthogonality relations the required formula was established. In [2] Faraut considered a general unimodular
group $ G $ and proved a proposition from which Gutzmer's formula can be deduced for noncompact Riemannian symmetric spaces [3] and Heisenberg groups [15].

Thus in all the previous settings the basic functions appearing in the Fourier
series or transform are matrix coefficients of certain irreducible unitary representations of the underlying group. Contrary to this, the Hermite functions do not occur as matrix coefficients. However, the Hermite functions are used to calculate the matrix coefficients associated to Schr\"odinger representations of
$ \H^n $ resulting in special Hermite or Laguerre functions. This explains why
the representation $ \pi(z,w) $ occurs in our Gutzmer's formula. The close relationship between Hermite and Laguerre functions are then used to derive the
Gutzmer's formula.

\section{\bf Preliminaries}
\setcounter{equation}{0}

In this section we collect some relevant information  about special
Hermite functions and prove some results that are required in the next
section. We closely follow the notations used in [12] and
[13] and we refer the reader to these monographs for more details.

Let $ \Phi_\alpha, \alpha \in \N^n $ be the Hermite functions on $ \R^n $
normalised so that their $ L^2 $ norms are one. These are eigenfunctions of 
the Hermite operator $ H = -\Delta +|x|^2 $ with eigenvalues $ (2|\alpha|+n).$
On finite linear combinations of such functions we can define certain operators $ \pi(z,w) $ where $ z, w \in ~\C^n $ as follows:
$$ \pi(z,w)\Phi_\alpha(\xi) = e^{i(z \cdot \xi+\frac{1}{2}z \cdot w)}
\Phi_\alpha
(\xi+w) $$ where $ z \cdot \xi = \sum_{j=1}^n z_j \xi_j $ and 
$ z \cdot w =     
\sum_{j=1}^n z_j w_j $. Note that $ \Phi_\alpha(\xi)= H_\alpha(\xi)
e^{-\frac{1}{2}|\xi|^2} $ where $ H_\alpha $ is a polynomial on $ \R^n $ and 
for $ z \in \C^n $ we define $ \Phi_\alpha(z) $ to be $H_\alpha(z) 
e^{-\frac{1}{2}z^2} $ where $ z^2 = z \cdot z.$ The special Hermite functions
$ \Phi_{\alpha,\beta}(z,w) $ are then defined by
$$
 \Phi_{\alpha,\beta}(z,w) = (2\pi)^{-\frac{n}{2}}( \pi(z,w)\Phi_\alpha,
\Phi_\beta).
$$
The restrictions of $ \Phi_{\alpha,\beta}(z,w) $ to $ \R^n\times\R^n $ are
usually called the special Hermite functions and the family
$\{\Phi_{\alpha,\beta}(x,u): \alpha, \beta \in \N^n \} $ forms an orthonormal
basis for $ L^2(\C^n).$

As we have mentioned in the introduction the operators $ \pi(z,w) $ are related
to the Schr\"odinger representation $ \pi_1 $ of the Heisenberg group $ \H^n.$
Recall that $ \H^n = \R^n \times \R^n \times \R $ is equipped with the group
law $ (x,u,t)(x',u',t') = (x+x',u+u',t+t'+\frac{1}{2}(u\cdot x'-x\cdot u')).$
For each  nonzero real number $ \lambda $ we have a representation of $ \H^n $
realised on $ L^2(\R^n) $ given by
$$ \pi_\lambda(x,u,t)\varphi(\xi) = e^{i\lambda t} e^{i\lambda(x\cdot \xi+\frac
{1}{2}x\cdot u)}\varphi(\xi+u).$$ Thus $ \pi(x,u) = \pi_1(x,u,0) $ and it defines a projective representation of $ \R^n \times \R^n.$

For $ (z,w) \in \C^{2n} $ the operators $ \pi(z,w) $ are not even bounded on
$ L^2(\R^n). $ However, they are densely defined and satisfy
$$
   \pi(z,w)\pi(z',w') = e^{\frac{i}{2}(z' \cdot w-z \cdot w')}\pi(z+z',w+w').
$$
Moreover,
$$ (\pi(iy,iv)\Phi_\alpha, \Phi_\beta) = (\Phi_\alpha,\pi(iy,iv)\Phi_\beta).$$
This means that $ \pi(iy,iv) $ are self adjoint operators. We need to 
calculate the $ L^2 $ norms of $ \pi(z,w)\Phi_\alpha.$ Let $ L_k^{n-1}$ be 
Laguerre polynomials of type $ (n-1) $ and define the Laguerre functions  
$ \varphi_k $ by
$$ \varphi_k(x,u) = L_k^{n-1}(\frac{1}{2}(x^2+u^2))e^{-\frac{1}{4}
(x^2+u^2)}.$$
Then it is known that
$$ \varphi_k(x,u) = (2\pi)^{n/2}\sum_{|\alpha|=k} \Phi_{\alpha,\alpha}(x,u).$$
These functions have a natural holomorphic extension to 
$ \C^n \times \C^n $ denoted by the same symbol:
$$
 \varphi_k(z,w) = (2\pi)^{n/2}\sum_{|\alpha|=k} \Phi_{\alpha,\alpha}(z,w).
$$

\begin{lem} For any $ z = x+iy, w = u+iv \in \C^n $ and $ \alpha \in \N^n $ we 
have
$$ \int_{\R^n} |\pi(z,w)\Phi_\alpha(\xi)|^2 d\xi = (2\pi)^{\frac{n}{2}}
e^{(u\cdot y -v\cdot x)} \Phi_{\alpha,\alpha}(2iy,2iv).$$
\end{lem}
{\bf Proof:} It is enough to prove the result in one dimension. Recall 
Mehler's formula satisfied by the Hermite functions $ h_k $ on $ \R $:
$$ \sum_{k=0}^\infty h_k(\xi)h_k(\eta) r^k = \pi^{-\frac{1}{2}}
e^{-\frac{1}{2}\frac{1+r^2}{1-r^2}(\xi^2+\eta^2)+\frac{2r}{1-r^2}\xi \eta}$$
valid for all $ r $ with $ |r| <1.$ The formula is clearly valid even if 
$ \xi $ and $ \eta $ are complex. A simple calculation shows that 
$$ \sum_{k=0}^\infty r^k |\pi(z,w)h_k(\xi)|^2 $$
$$ = \pi^{-\frac{1}{2}}(1-r^2)^{-\frac{1}{2}} e^{-(uy+vx)}e^{\frac{1+r}{1-r}
v^2}e^{-\frac{1-r}{1+r}(\xi+u)^2}e^{-2y\xi}.$$
Integrating both sides with respect to $ \xi $ we obtain
$$  \sum_{k=0}^\infty r^k \int_{\R} |\pi(z,w)h_k(\xi)|^2  d\xi $$
$$= (1-r)^{-1} e^{(uy-vx)} e^{\frac{1+r}{1-r}(y^2+v^2)}.$$
We now recall that the generating function for the Laguerre functions $ \varphi_k(x,u) $ when $ n =1 $ reads as
$$ \sum_{k=0}^\infty r^k \varphi_k(x,u) = (1-r)^{-1} e^{-\frac{1}{4}
(x^2+u^2)}.$$ A comparison with this shows that
$$ \int_\R |\pi(z,w)h_k(\xi)|^2  d\xi =  e^{(uy-vx)} \varphi_k(2iy,2iv).$$ 
Since $ \Phi_{k,k}(x,u) = (2\pi)^{-\frac{1}{2}}\varphi_k(x,u) $ this proves 
the Lemma.

In the above lemma we have calculated the $ L^2 $ norm of 
$ \pi(z,w)\Phi_\alpha $ by integrating the generating function. We can also
calculate the norm by expanding $ \pi(z,w)\Phi_\alpha $ in terms of the
Hermite basis and appealing to the Plancherel theorem for Hermite expansions.
This leads to the following identity which is crucial for our main result.

\begin{lem} For any $ \alpha \in \N^n, z = x+iy, w = u+iv \in \C^n $ we have
$$ \sum_{\beta \in \N^n} |\Phi_{\alpha,\beta}(z,w)|^2 = (2\pi)^{\frac{-n}{2}} 
e^{(u\cdot y -v\cdot x)} \Phi_{\alpha,\alpha}(2iy,2iv).$$
\end{lem}
{\bf Proof:} We just have to recall that $ (\pi(z,w)\Phi_\alpha,\Phi_\beta)
= (2\pi)^{\frac{n}{2}}\Phi_{\alpha,\beta}.$

We also need some estimates on the holomorphically extended Hermite functions
on $ \C^n.$ Let us define $ \Phi_k(x,u) = \sum_{|\alpha| = k}\Phi_\alpha(x)
\Phi_\alpha(u) $ which is the kernel of the projection $ P_k.$ Note that $ 
\Phi_k $ extends to $ \C^n \times \C^n $ as an entire function. Using Mehler's
formula for Hermite functions and the generating function for Laguerre 
functions we can get the following representation of $ \Phi_k $ in terms
of Laguerre functions of type $ (n/2-1).$

\begin{lem}
$$ \Phi_k(z,w) = \pi^{-\frac{n}{2}}\sum_{j=0}^k (-1)^j L_j^{n/2-1}(\frac{1}{2}
(z+w)^2)L_{k-j}^{n/2-1}(\frac{1}{2}(z-w)^2) e^{-\frac{1}{2}(z^2+w^2)}$$
where $ z^2 = \sum_{j=1}^n z_j^2 $ and $ w^2 = \sum_{j=1}^n w_j^2.$
\end{lem}
{\bf Proof:} The Laguerre functions of type $ (n/2-1) $ are given by the
generating function
$$ \sum_{k=0} r^k L_k^{n/2-1}(\frac{1}{2}z^2)e^{-\frac{1}{4}z^2}
= (1-r)^{-n/2}e^{-\frac{1}{4}\frac{1+r}{1-r}z^2}.$$
A simple calculation shows that
$$ (1-r)^{-n/2}e^{-\frac{1}{4}\frac{1+r}{1-r}(z+w)^2}(1+r)^{-n/2}
e^{-\frac{1}{4}\frac{1-r}{1+r}(z-w)^2}$$ 
$$ = (1-r^2)^{-n/2} e^{-\frac{1}{2}\frac{1+r^2}{1-r^2}(z^2+w^2)+ 
\frac{2r}{1-r^2}zw}.$$
Comparing this with Mehler's formula and rewriting the left hand side as
a power series in $ r $ and then equating coefficients of $ r^k $ we 
obtain the lemma.

The above lemma has been already used by us in the study of Bochner-Riesz means
for multiple Hermite expansions. Here we need the above in order to get the
following estimate on $ \Phi_k(z,w).$

\begin{lem} For all $ z = x+iy \in \C^n $ and $ k = 1,2,... $ we have
$$ |\Phi_k(z,\bar{z})| \leq C(y) k^{\frac{3}{4}(n-1)}
 e^{2(k)^{\frac{1}{2}}|y|} $$ where $ C(y) $ is locally bounded.
\end{lem}
{\bf Proof:} From the previous lemma we have
$$ \Phi_k(z,\bar{z}) = \pi^{-\frac{n}{2}}\sum_{j=0}^k (-1)^j L_j^{n/2-1}
(2|x|^2)e^{-|x|^2} L_{k-j}^{n/2-1}(-2|y|^2) e^{|y|^2}.$$
We now make use of the following estimates on Laguerre functions. First of all
we know that
$$ |L_j^{n/2-1}(2|x|^2)e^{-|x|^2}| \leq C j^{n/2-1} $$ uniformly in $ x.$ On the other hand Perron's formula for Laguerre polynomials in the complex domain
(see Theorem 8.22.3 in Szego [11] ) gives us
$$ L_{j}^{n/2-1}(-2|y|^2) e^{|y|^2} \leq C(y) j^{\frac{(n-3)}{4}}
 e^{2(j)^{\frac{1}{2}}|y|} $$ valid for all $ |y| \geq 1.$ Since $ 
L_{j}^{n/2-1}(-2|y|^2) \leq L_{j}^{n/2-1}(-2) $ we have the same estimate for 
all values of $ y.$ These two estimates give the required bound on $ \Phi_k(z,
\bar{z}).$

We conclude the preliminaries with establishing some more notation. Let 
$ Sp(n,\R) $ stand for the symplectic group consisting of $ 2n \times 2n $ real
matrices that preserve the symplectic form $ [(x,u),(y,v)] = (u\cdot y-v\cdot x)
$ on $ \R^{2n} $ and have determinant one. Let $ O(2n,\R) $ be the orthogonal
group and we define $ K = Sp(n,\R)\cap O(2n,\R).$ Then there is a one to one
correspondence between $ K $ and the unitary group $ U(n) .$ Let $ \sigma 
= a+ib $ be an $ n \times n $ complex matrix with real and imaginary parts $ a $ and $ b.$ Then $ \sigma $ is unitary if and only if the matrix $ A = 
\begin{pmatrix}a & -b \cr b & a  
\end{pmatrix}$ is
in $ K.$ For these facts we refer to Folland [4]. By $ \sigma.(x,u) $ we denote the action of the correspoding matrix $ A $ on $ (x,u).$ This action has a natural extension to $ \C^n \times \C^n $ denoted by $ \sigma.(z,w) $ and is 
given by $ \sigma.(z,w) = (a.z-b.w, a.w+b.z) $ where $ \sigma = a+ib.$ 
For example, when $ n = 1 $ and $ \sigma = e^{i\theta} $ we see that the 
corresponding matrix $ A $ is $
\begin{pmatrix} \cos\theta & -\sin\theta \cr
\sin\theta & \cos\theta .
\end{pmatrix}$ Given $ \theta = (\theta_1,....,\theta_n) \in
\R^n $ we denote by $ k(\theta) $ the diagonal matrix in $ U(n) $ with entries
$ e^{i\theta_j}.$ We denote by $ d\sigma $ the normalised Haar measure on
$ K $ and by $ d\theta $ the Lebesgue measure 
$ d\theta_1 d\theta_2....d\theta_n.$

\section{\bf The main results}
\setcounter{equation}{0}

Having set up notation and collected relevant results on special Hermite functions we are now ready to prove our main results. We begin with

\begin{thm} Let $ f \in L^2(\R^n) $ be such that $ \|P_kf\|_2 \leq C_t 
e^{-2k^{\frac{1}{2}}t} $ for all $ t > 0 $ and $ k \in \N.$ Then $ f $ has a
holomorphic extension $ F $ to $ \C^n $ and we have the following formula
for any  $ z = x+iy, w= u+iv \in \C^n $:
$$ \int_{\R^n} \int_K |\pi(\sigma.(z,w))F(\xi)|^2 d\sigma d\xi $$
$$ =  e^{(u\cdot y - v \cdot x)} \sum_{k=0}^\infty \frac{k!(n-1)!}{(k+n-1)!}
\varphi_k(2iy,2iv) \|P_kf\|_2^2 .$$
\end{thm}
{\bf Proof:} Consider the Hermite expansion of the function $ f $ given by
$$  f(x) = \sum_{k=0}^\infty \sum_{|\alpha|=k}(f,\Phi_\alpha)\Phi_\alpha(x).$$
By Cauchy-Schwarz inequality
$$ |\sum_{|\alpha|=k}(f,\Phi_\alpha)\Phi_\alpha(x+iy)|^2 \leq 
\Phi_k(x+iy,x-iy) \|P_kf\|_2^2 .$$
In view of Lemma 2.4 the hypothesis on $ f $ allows us to
conclude that the series
$$  \sum_{k=0}^\infty \sum_{|\alpha|=k}(f,\Phi_\alpha)\Phi_\alpha(x+iy) $$
converges uniformly over compact subsets of $ \C^n $ and hence $ f $ extends 
to an entire function $ F $ on $ \C^n.$

Let $ D $ be the subgroup of $ K $ consisting of $ 2n \times 2n $ matrices 
associated to the elements $ k(\theta) \in U(n).$ We claim that it is enough
to prove
$$
 (2\pi)^{-n} \int_{\R^n} \int_D |\pi(k(\theta).(z,w))F(\xi)|^2 d\theta d\xi $$
$$ = (2\pi)^{n/2} e^{(u\cdot y - v \cdot x)} \sum_{\alpha \in \N^n}
 \Phi_{\alpha,\alpha}(2iy,2iv)|(f,\Phi_\alpha)|^2.
$$
To see the claim, suppose we have the above formula. Then writing
$$ \int_{\R^n} \int_K |\pi(\sigma.(z,w))F(\xi)|^2 d\sigma d\xi $$
$$  =  (2\pi)^{-n}\int_{\R^n} \int_D \int_K |\pi(k(\theta)\sigma.(z,w))
F(\xi)|^2 d\sigma d\theta d\xi $$
we get
$$ \int_{\R^n} \int_K |\pi(\sigma.(z,w))F(\xi)|^2 d\sigma d\xi $$
$$ = (2\pi)^{n/2} e^{(u'\cdot y' - v' \cdot x')} \sum_{\alpha \in \N^n}
 \Phi_{\alpha,\alpha}(2iy',2iv')|(f,\Phi_\alpha)|^2 $$
where $ (z',w') = \sigma.(z,w).$ Since the action of $ \sigma $ preserves the
symplectic form we have $ e^{(u\cdot y - v \cdot x)} = 
e^{(u'\cdot y' - v' \cdot x')}.$ Thus we are left with proving
$$ (2\pi)^{n/2} \int_K \Phi_{\alpha,\alpha}(\sigma.(2iy,2iv)) d\sigma = 
\frac{k!(n-1)!}{(k+n-1)!}\varphi_k(2iy,2iv) $$ 
whenever $ |\alpha| =k.$ But this is a well known fact. A representation 
theoretic proof of this can be found in Ratnakumar et al [10]. 

(Another way to see this is the following. The functions
$ \Phi_{\alpha,\alpha}(x,u) $ are eigenfunctions of the special Hermite 
operator $ L $ with eigenvalue $ (2|\alpha|+n).$ And hence the function
$ \int_K \Phi_{\alpha,\alpha}(\sigma.(x,u)) d\sigma $ is a radial eigenfunction
of the same operator. But any bounded radial eigenfunction with eigenvalue 
$(2k+n) $ is a constant 
multiple of $ \varphi_k(x,u).$ This proves that
$$ (2\pi)^{n/2} \int_K \Phi_{\alpha,\alpha}(\sigma.(x,u)) d\sigma = 
\frac{k!(n-1)!}{(k+n-1)!}\varphi_k(x,u) $$ and hence they are same on
$ \C^n \times \C^n $ as well.)

We now turn our attention to prove the formula for the action of $ D.$ The 
idea is to expand the operator valued function $ \pi(k(\theta).(z,w)) $ into
a Fourier series. Defining
$$ \pi_m(z,w)F(\xi) = (2\pi)^{-n} \int_D \pi(k(\theta).(z,w))F(\xi)
e^{-im\cdot \theta} d\theta $$ we have the expansion
$$ \pi(k(\theta).(z,w))F(\xi) = \sum_{m \in \Z^n} 
\pi_m(z,w)F(\xi) e^{im\cdot \theta}.$$ By the orthogonality of the Fourier 
series we obtain
$$ (2\pi)^{-n} \int_{\R^n}\int_D |\pi(k(\theta).(z,w))F(\xi)|^2 d\theta d\xi $$
$$ = \sum_{m \in \Z^n} \int_{\R^n}|\pi_m(z,w)F(\xi)|^2 d\xi.$$ In calculating 
the $ L^2 $ norm of $ \pi_m(z,w)F $ we make use of another property of special
Hermite functions, namely that $ \Phi_{\alpha,\beta}(x,u) $ is $ 
(\beta-\alpha)-$ homogeneous. By this we mean
$$ \Phi_{\alpha,\beta}(k(\theta).(x,u)) = e^{i(\beta-\alpha)\cdot \theta}
\Phi_{\alpha,\beta}(x,u).$$ A proof of this can be found in [12] (see 
Proposition 1.4.2).

Expanding $ f $ in terms of the Hermite basis we see that
$$ \pi_m(z,w)F  =  \sum_{\alpha,\beta} (f,\Phi_\alpha)
(\pi_m(z,w)\Phi_\alpha,\Phi_\beta) \Phi_\beta.$$ 
But
$$ (\pi_m(x,u)\Phi_\alpha,\Phi_\beta) = (2\pi)^{-n/2}\int_D 
\Phi_{\alpha,\beta}(k(\theta).(x,u))e^{-im\cdot \theta}d\theta = 0 $$
unless $ \beta = \alpha +m $ due to the homogeneity properties of the special
Hermite functions. Therefore, the expansion of $ \pi_m(z,w)F $ reduces to
$$ \pi_m(z,w)F  = (2\pi)^{n/2} \sum_{\alpha \in \N^n} (f,\Phi_\alpha)
\Phi_{\alpha,\alpha +m}(z,w) \Phi_{\alpha +m}.$$ This leads us to
$$ \|\pi_m(z,w)F\|_2^2 = (2\pi)^n \sum_{\alpha \in \N^n} |(f,\Phi_\alpha)|^2
|\Phi_{\alpha,\alpha +m}(z,w)|^2.$$
Thus we have proved
$$ (2\pi)^{-n} \int_{\R^n}\int_D |\pi(k(\theta).(z,w))F(\xi)|^2 d\theta d\xi $$
$$ = (2\pi)^n \sum_{m \in \Z^n} \sum_{\alpha \in \N^n} |(f,\Phi_\alpha)|^2
|\Phi_{\alpha,\alpha +m}(z,w)|^2. $$
This proves our claim since the sum over $ m \in \Z^n $ is precisely
$ (2\pi)^{-n/2} e^{(u\cdot y - v \cdot x)} \Phi_{\alpha,\alpha}(2iy,2iv)$ 
in view of Lemma 2.2. Hence the proof of the theorem is complete.

The above theorem has a natural converse which we state and prove now. Together
they prove Theorem 1.1 stated in the introduction. In the proof of the above
theorem the hypothesis on the Hermite projections of $ f $ are used twice. 
First we used the estimates to conclude that $ f $ has an entire extension to
$ \C^n.$ Then we used them to show that the sum and the integral appearing
in the above theorem are finite. In the next theorem we begin with an entire
function for which the integral is finite and obtain the estimates on the
projections.

\begin{thm} Let $ F $ be an entire function on $ \R^n $ for which the integral
$$ \int_{\R^n} \int_K |\pi(\sigma.(z,w))F(\xi)|^2 d\sigma d\xi $$
is finite for all $ z, w \in \C^n.$ Then $ \|P_kf\|_2 \leq C_t
 e^{-2k^{\frac{1}{2}}t} $ for all $ t > 0.$
\end{thm}
{\bf Proof:} We proceed as in the proof of the previous theorem. Since $ F $ is
holomorphic $ \pi(z,w)F $ makes sense. As before, for almost every $ \sigma \in
U(n) $ we have
$$ \int_{\R^n} \int_D |\pi(k(\theta)\sigma.(z,w))F(\xi)|^2 d\theta d\xi 
< \infty .$$
Expanding the operator $ \pi(k(\theta).(z,w)) $ into Fourier series and proceeding exactly as in the previous theorem and noting that at each stage the
resulting sums are finite we get the Gutzmer's formula, namely the integral
in the theorem is equal to
$$ e^{(u\cdot y-v\cdot x)} \sum_{k=0}^\infty \frac{k!(n-1)!}{(k+n-1)!}
\varphi_k(2iy,2iv)\|P_kf\|_2^2 $$
and hence the sum is finite. Now Perron's formula for Laguerre functions
on the negative real axis also gives lower bounds. That is to say, the
Laguerre functions $ \varphi_k(2iy,2iv) $ behave like $ e^{2(k)^{\frac{1}{2}}
(|y|^2+|v|^2)^{\frac{1}{2}}}.$ In view of this we immediately get the decay
estimates on the projections $ P_kf.$

\section{\bf Some consequences}
\setcounter{equation}{0}

In this section we deduce some interesting consequences of our Gutzmer's
formula. First we obtain the characterisation of the image of $ L^2(\R^n) $ 
under the Hermite semigroup mentioned in Corollary 1.2. As we have pointed 
out earlier the result is not new but we give a different proof.

Consider the heat kernel $ p_t(y,v) $ associated to the special Hermite 
operator which is explicitly given by
$$ p_t(y,v) = (2\pi)^{-n} (\sinh(t))^{-n} 
e^{-\frac{1}{4}\coth(t)(|y|^2+|v|^2)}.$$
We now look at the integral
$$ \int_{\R^n} \left( \int_{\R^{2n}} |\pi(iy,iv)f(\xi)|^2 p_{2t}(2y,2v)dydv 
\right) d\xi.$$
Since the function $ p_t(y,v) $ and the Lebesgue measure $ dy dv $ are both
invariant under the action of the group $ K $ we can rewrite the above
integral as
$$ \int_{\R^{2n}}\left( \int_{\R^n} \int_K |\pi(\sigma.(iy,iv))f(\xi)|^2 
d\sigma d\xi \right) p_{2t}(2y,2v)dydv .$$
In view of Gutzmer's formula the above reduces to
$$ \sum_{k=0}^\infty \frac{k!(n-1)!}{(k+n-1)!} \left(\int_{\R^{2n}}
\varphi_k(2iy,2iv)p_{2t}(2y,2v) dydv \right) \|P_kf\|_2^2.$$
We now make use of the fact that
$$ \frac{k!(n-1)!}{(k+n-1)!}\int_{\R^{2n}}
\varphi_k(2iy,2iv)p_{2t}(2y,2v) dydv = e^{2(2k+n)t} $$ which we have 
established in [15] (see Lemma 6.3).

Therefore, replacing $ f $ by $ e^{-tH}f $ we have established
$$ \int_{\R^{2n}} \left( \int_{\R^n} |\pi(iy,iv)e^{-tH}f(\xi)|^2 d\xi \right)
p_{2t}(2y,2v) dy dv $$
$$ = \sum_{k=0}^\infty \|P_kf\|_2^2 = \int_{\R^n} |f(\xi)|^2 d\xi .$$
Writing $ F $ for $ e^{-tH}f $ a simple calculation shows that the above 
integral is equal to
$$ (2\pi \sinh(2t))^{-n} \int_{\R^{2n}} \left( \int_{\R^n} |F(\xi+iv)|^2 e^{-2y\cdot \xi} e^{-\coth(2t)(|y|^2+|v|^2)} dy \right) d\xi dv.$$
Performing the integration with respect to $ y $ we see that the above is
nothing but
$$ \int_{\R^{2n}}|F(\xi+iv)|^2 U_t(\xi,v) d\xi dv .$$ This completes the proof of Corollary 1.2.

We remark that if we have only assumed the estimate $ \|P_kf\|_2 \leq C 
e^{-2k^{\frac{1}{2}}t} $  for some $ t > 0 $ ( not for all $ t $ as in
Theorem 3.1) then the proof of Theorem 3.1 shows that $ f $ can be extended
as a holomorphic function to cetain tube domain $ \Omega_t = \{ z \in \C^n:
|y| < t \} $ and still we have Gutzmer's formula as long as $ |y|^2+|v|^2 < t^2.$ We may think of Gutzmer's formula as a characterisation of the image of
$ L^2(\R^n) $ under the Hermite-Poisson semigroup $ e^{-tH^{\frac{1}{2}}}.$ 
Compare this with the results of Janssen and Eijndhoven [5] on the growth of 
Hermite coefficients.

Another interesting consequence of the Gutzmer's formula is the following orthogonality relations for Hermite functions on $ \C^n.$ Polarising Gutzmer we obtain
$$ \int_{\R^n} \int_K \pi(\sigma.(z,w))F(\xi) 
\overline{\pi(\sigma.(z,w))G(\xi)} d\sigma d\xi $$
$$ =  e^{(u\cdot y - v \cdot x)} \sum_{k=0}^\infty \frac{k!(n-1)!}{(k+n-1)!}
\varphi_k(2iy,2iv) (P_kf,P_kg) .$$
Specialising to Hermite functions we get the following result which, to our knowledge, seems to be new.

\begin{cor} For any $ z,w \in \C^n $ and $ \alpha,\beta \in \N^n$ we have
$$\int_{\R^n} \int_K \pi(\sigma.(z,w))\Phi_\alpha(\xi)
\overline{\pi(\sigma.(z,w))\Phi_\beta(\xi)} d\sigma d\xi $$
$$ = e^{(u\cdot y - v \cdot x)}\frac{k!(n-1)!}{(k+n-1)!}\varphi_k(2iy,2iv) 
\delta_{\alpha,\beta}.$$
\end{cor}

The above shows that in the one dimensional case the Hermite functions $ h_k $
satisfy the following relations. The choice $ z = i\eta, w = 0 $ gives
$$ \int_{\R} \int_0^{2\pi} e^{-2\xi \eta \cos \theta} 
h_k(\xi+i\eta \sin \theta)\overline{h_j(\xi+i\eta \sin \theta)} d\theta d\xi $$
$$ = (2\pi) L_k^0(-2\eta^2)e^{\eta^2} \delta_{k,j}.$$
The choice $ z = \eta, w = i\eta $ leads to
$$ \int_{\R} \int_0^{2\pi}e^{2\xi \eta \sin \theta - \eta^2 \cos(2\theta)}
h_k(\xi+i\eta e^{-i\theta})\overline{h_j(\xi+i\eta e^{-i\theta})} 
d\theta d\xi $$
$$ = (2\pi) L_k^0(-2\eta^2) \delta_{k,j}.$$
Other interesting relations in higher dimensional cases can be obtained by 
suitable choices of $ z, w $ and also by choosing various subgroups of $ K.$


\begin{thebibliography}{99}

\bibitem{By} D-W. Byun, \textit{Inversions of Hermite semigroup}, Proc.
A.M.S. {\bf 118}(1993), 437-445.


\bibitem{} J. Faraut, \textit{Formule de Gutzmer pour la complexification 
d'un espace Riemannien symetrique}, Rend. Mat. Acc. Lincei s.9, {\bf v.13}, 
233-241 (2002).

\bibitem{} J. Faraut, \textit{Analysis on the crown of a Riemannian symmetric
space},  Amer. Math. Soc. Transl. {\bf 210}, no. 2, 99-110 (2003).

\bibitem{} G. B. Folland, \textit{ Harmonic analysis in phase space}, Ann.
Math. Stud. {\bf 122}, Princeton Univ. Press, Princeton (1989).

\bibitem{} A. J. E. M. Janssen and S. J. L. van Eijndhoven, \textit{ Spaces 
of type W, growth of Hermite coefficients, Wigner distribution, and Bargmann 
transform}, J. Math. Anal. Appli.{\bf 152}, 368-390 (1990).


\bibitem{} D. Karp, \textit{ Square summability with geometric weight 
for classical orthogonal expansions}, Advances in Analysis, Ed. H. G. W. Begehr
et al, World Scientific, 407-421 (2005).

\bibitem{} B. Kr\"otz, S. Thangavelu and Y. Xu, \textit{ The heat kernel 
transform for the Heisenberg group}, J. Funct. Anal. {\bf 225}, no.2, 301-336 
(2005).

\bibitem{} B. Kr\"otz, G. Olafsson and R. Stanton, \textit{ The image of the 
heat kernel transform on Riemannian symmetric spaces of noncompact type}, Int.
Math. Res. Notes, no. {\bf 22}, 1307-1329 (2005).

\bibitem{} M. Lassalle, \textit{Series de Laurent des fonctions holomorphes 
dans la
complexification d'un espace symetrique compact}, Ann. Sci. Ecole Norm. Sup.
{\bf 11} (1978),167-210.

\bibitem{} P. K. Ratnakumar, R. Rawat and S. Thangavelu, \textit{ A 
restriction theorem for the Heisenberg motion group}, 
Stud. Math. {\bf 126(1)}, 1-12 (1997).


\bibitem{} G. Szeg\"o, \textit{Orthogonal polynomials}, Amer. Math. Soc. 
Colloq. Publi., Providence, RI (1967).


\bibitem{} S. Thangavelu, \textit{Harmonic analysis on the Heisenberg group},
Prog. in Math. Vol. {\bf 159}, Birkh\"auser, Boston (1998).

\bibitem{} S. Thangavelu, \textit{ An introduction to the uncertainty 
principle},Prog. in Math. Vol. {\bf 217}, Birkh\"auser, Boston (2004).

\bibitem{} S. Thangavelu, \textit{ Hermite and Laguerre semigroups: some recent
developments}, CIMPA- VENEZUELA Lecture Notes (2006).

\bibitem{} S. Thangavelu, \textit{Gutzmer's formula and Poisson integrals on the
Heisenberg group}, Pacific J. Math. (to appear).

\end{thebibliography}
\end{document}